\documentclass [18pt]{article}
\usepackage {graphicx}

\begin{document}
\textbf{Application of Quantum Theory to Super-parametric Density
Estimation}

Yeong-Shyeong Tsai

Department of Applied Mathematics, National Chung Hsing
University, Taichung,Taiwan

\textbf{Abstract}

Since the consistency of maximum likelihood estimator has been proved, the
only problem which is left is the problem of optimization. In last century,
it was found that some splines were very useful. From Stone-Weirstrass
theorem, we can approximate continuous functions by the polynomials and
hence we can construct the estimator by using the splines. Therefore, it
might not be so important to stress the difference between the parametric
approach and nonparametric approach. Usually, a nonlinear optimization
problem is not so easy to solve and it is assumed that the optimization
problem can be solved by existent packages. From the view point of
mathematics, the results of the optimization problem should be verified or
reinvestigated because the nonlinear optimization problem is not simple. It
seems that the nonlinear optimization play an important role in density
estimation. Though nonlinear equations must be solved in most optimization
problems, we will show how a optimization problem can be solved by finding
the solution of systems of linear equations. Basing on this approach, the
optimization problem can be solved by solving a quadratic equation finally.
Some numerical examples are studied as well. From the figures, it can be
found this is a good approach on density function estimation.

\textbf{Introduction}

The problem of density estimation is to estimate the density function $f$ by
a set of observations, $x_1 $, $x_2 $,\ldots , $x_m $. Roughly, the function
with parameters, denoted by the notation$f(x,\theta )$, is called the
estimator of $f$. We assume that there is a family of functions, say

\begin{equation}
\label{eq1}
\Im = \{f(x,\theta );\theta \in R^n\},
\end{equation}

\noindent
and $f \in \Im $. The likelihood function $l$ is defined

\begin{equation}
\label{eq2}
l = \prod\limits_{j = 1}^m {f(x_j ,\theta )} .
\end{equation}

From the work of the statistician, the information of $f$ can be obtained by
maximizing the likelihood function of the density estimator [1]. Usually, it
is not a simple work to solve the nonlinear equations. So far, we know how
to solve a single linear equation, a system of linear equations and a single
quadratic equation. In this paper, the optimization problem is transformed
to system of linear equations first. Basing on the approach, the
optimization problem is transformed to a single quadratic equation. Finally,
we can solve the optimization problem effectively. The work of
transformation is not so simple though the idea is simple. Besides, the
undesired roughness of nonparametric estimator is a serious problem. Since
our approach is expected to estimate the density function of general cases,
this serious problem must be studied in the same time. In the middle of
20$^{th}$ century, several splines were studied. There are many applications
of these splines such as computer aid design of cars [2], curve fitting in
statistics, computation of energy levels of multi-electron atoms etc [3].
These splines were introduced to diminish the oscillations of the curve
which is obtained by the method of traditional polynomial curve fitting.
Therefore, these splines can remove the roughness of density estimators. It
is possible to solve these problems in the same time. Anyone who knows the
elementary calculus [4] or second year calculus [5] is able to understand
this paper.

\textbf{Parzen windows}

In order to avoid the difficulty of the nonlinear optimization problem, the
orthogonal polynomials are used in most nonparametric methods. It seems that
the orthogonal function will introduce more roughness. In order to avoid
introducing the roughness, the orthonormal basis is abandoned and the Parzen
window functions [6], nonnegative functions, are adopted. Let $\delta $ be
the Dirac delta function. The Dirac delta function is a generalized
function,

\begin{equation}
\label{eq3}
\delta (x) = 0
\end{equation}

\noindent
when $x \ne 0$, and

\begin{equation}
\label{eq4}
\int_{ - \infty }^\infty {\delta (x) = 1} .
\end{equation}

If it is necessary, then we shall consider the Dirac delta function as a
linear functional defined on a function space [7]. Intuitively, we can start
from the following identity

\begin{equation}
\label{eq5}
f(x) = \int {\delta (x - t)f(t)dt} .
\end{equation}

Here, $f$ is the probability density function which will be estimated by
observations $x_1 $, $x_2 $, $x_3 $,\ldots , $x_m $. Let $\widetilde{f}$ be
the estimator of $f$. If the integration of (\ref{eq5}) can be approximated by
summation, then the estimator is

\begin{equation}
\label{eq6}
\widetilde{f}(x) = \sum\limits_{i = 1}^n {c_i } \varphi _i (x),
\end{equation}

\begin{equation}
\label{eq7}
\varphi _i (x) \ge 0.
\end{equation}

Usually, $\varphi _i $ are called window functions or kernel functions. It
seems that (\ref{eq1})-(\ref{eq5}) can be ignored. We can start from the estimator which is
defined in (\ref{eq6}). If we can determine value of $c_i $ properly, then the
estimator is obtained. Though there are many window functions which are
available [2] [3], we find that Bernstein polynomial is good a candidate.
Clearly, it must be that

\begin{equation}
\label{eq8}
\int {\widetilde{f}(x)dx = 1} .
\end{equation}

Let

\begin{equation}
\label{eq9}
p_i = \int {\varphi _i } (x)dx.
\end{equation}

Let $l$ be defined

\begin{equation}
\label{eq10}
l = \prod\limits_{j = 1}^m {\widetilde{f}(x_j )} .
\end{equation}

Here, $l$ is the likelihood function. From the works of statisticians, the
value of $c_i $ can be determined by maximizing the likelihood function [1].
The problem is to maximize $l$ subjected to the constraints,

\begin{equation}
\label{eq11}
\sum\limits_{i = 1}^n {p_i c_i = 1} ,
\end{equation}

\noindent
and

\begin{equation}
\label{eq12}
0 \le c_i ,
\quad
i = 1,2,...,n.
\end{equation}

Mathematically, since $c_i $ are going to be determined, if we redefine
$\widetilde{f}$,

\begin{equation}
\label{eq13}
\widetilde{f}(x) = \sum\limits_{i = 1}^n {(c_i / p_i )} \varphi _i (x),
\end{equation}

Then the constraints become

\begin{equation}
\label{eq14}
\sum\limits_{i = 1}^n {c_i = 1}
\end{equation}

\noindent
and

\begin{equation}
\label{eq15}
0 \le c_i ,
\quad
i = 1,2,...,n.
\end{equation}

Generally, this problem should be solved by Kuhn-Tucker Theorem [8]. Like
many mathematical theorems, both Kuhn-Tucker and Lagrange theories are not
constructive and the non-constructive results can be traced back to the last
axiom of real number, axiom of completeness [4]. In physics, the orthogonal
functions are very useful. In order to solve the nonlinear optimization of
density estimation, the orthogonal functions were adopted by nonparametric
approach. Quantum mechanics was discussed in the paper of Good and Gaskins
1971 [9] [10]. Since then most, if not all, nonparametric density estimators
were built on the orthogonal functions which were inferred from quantum
mechanics directly or indirectly. Let $V$ be a vector space over the field
of complex numbers. Let $v_1 $, $v_2 $,\ldots , $v_n $ be orthonormal basis
of $V$. Let $z$ be a unit vector of $V$. Roughly speaking, If $z =
\sum\limits_{i = 1}^n {c_i v_i } $, then $\left| {c_i } \right|^2$ is
interpreted as the probability that $z$ might be $v_i $ in quantum theory.
Hence $c_i $ is called the probability amplitude. In statistics, the real
numbers work well. Therefore, the complex field is replaced by the real
field. The difference between the probability density and the probability
amplitude is clear and simple. Mathematically or symbolically, the symbol
$c_i $ is replaced by $c_i^2 $. Obviously, the constraints become

\begin{equation}
\label{eq16}
\sum\limits_{i = 1}^n {c_i^2 = 1} .
\end{equation}

Only one constraint is left. Then the Lagrange's multiplier technique can be
applied easily. Clearly, optimization problem play an important role in
density function estimation. In order to avoid the difficulties, we will
follow the approach of quantum theory and the concept of the probability
amplitude is adopted. But we will use the result of Stone-Weirstrass theorem
instead of the orthogonal functions

\textbf{Optimization on the compact manifold}

Since the likelihood function $l$ is a $C^\infty $ function of $c_i $
defined on the compact subset of $R^n$, $l$ must has maximum value on the
sphere. Let observations be $x_1 $, $x_2 $, $x_3 $,\ldots , $x_m $.

Let

\begin{equation}
\label{eq17}
l_j = \sum\limits_{i = 1}^n {c_i^2 \varphi _i (x_j )} .
\end{equation}

Let

\begin{equation}
\label{eq18}
l = \prod\limits_{j = 1}^m {l_j } .
\end{equation}

Let

\begin{equation}
\label{eq19}
\sum\limits_{i = 1}^n {c_i c_i } = r
\end{equation}

\noindent
be the constraint. Now, we start to solve the optimization problem.

Let

\begin{equation}
\label{eq20}
a_{ij} = \varphi _i (x_j ).
\end{equation}

Let

\begin{equation}
\label{eq21}
l_L = l - \lambda (\sum\limits_{i = }^n {c_i c_i } - r).
\end{equation}

By the method of Lagrange's multiplier, we have

\begin{equation}
\label{eq22}
\frac{\partial l_L }{\partial c_k } = 0.
\end{equation}

From (\ref{eq17})- (\ref{eq22}), we get

\begin{equation}
\label{eq23}
l\sum\limits_{j = 1}^m {\frac{a_{kj} c_k }{l_j }} - \lambda c_k = 0,
\quad
k = 1,2,...,n.
\end{equation}

Multiplying (\ref{eq23}) by $c_k $ and taking the summation of index $k$, we get

\begin{equation}
\label{eq24}
\sum\limits_{k = 1}^n {(l\sum\limits_{j = 1}^m {\frac{a_{kj} c_k }{l_j }} -
\lambda c_k )c_k } = 0,
\end{equation}

Interchanging the summations, we get

\begin{equation}
\label{eq25}
\sum\limits_{j = 1}^m {(l\sum\limits_{k = 1}^n {\frac{a_{kj} c_k c_k }{l_j
}} - \lambda c_k c_k )} = 0.
\end{equation}

From (\ref{eq17}), (\ref{eq18}), (\ref{eq19}), (\ref{eq20}) and (\ref{eq25}), we get

\begin{equation}
\label{eq26}
lm - \lambda r = 0.
\end{equation}

Hence

\begin{equation}
\label{eq27}
l = \lambda r / m.
\end{equation}

Substituting (\ref{eq27}) into (\ref{eq23}), we get

\begin{equation}
\label{eq28} c_k (\sum\limits_{j = 1}^m {\frac{a_{kj} r}{ml_j }} -
1) = 0, \quad k = 1,2,...,n.
\end{equation}

Clearly, either

\begin{equation}
\label{eq29}
c_k = 0,
\end{equation}

\noindent
or

\begin{equation}
\label{eq30}
(\sum\limits_{j = 1}^m {\frac{a_{kj} r}{ml_j }} - 1) = 0,
\quad
k = 1,2,...,n.
\end{equation}

It should be emphasized that (\ref{eq29}) and (\ref{eq30}) are not mutually exclusive. In
order to linearization the equations (\ref{eq30}), we take some transformations of
variables.

Let

\begin{equation}
\label{eq31}
y_j = r / (ml_j ).
\end{equation}

Substituting (\ref{eq31}) into (\ref{eq30}), we get

\begin{equation}
\label{eq32}
\sum\limits_{j = 1}^m {a_{kj} y_j = 1}
\end{equation}

Multiplying both side (\ref{eq32}) by a constant $\theta ^2$, we have

\begin{equation}
\label{eq33}
\sum\limits_{j = 1}^m {a_{kj} \theta ^2y_j = \theta ^2} .
\end{equation}

Since the solutions of (\ref{eq29}) are not affected by any factor, the constant
$\theta $ is introduced in equations (\ref{eq32}) to fit the constraint. It seems
that the constant $\theta $ is a redundancy because that $\theta $ must be
1. Later, it will be found in following lemmas and theorems that $\theta $
play an important role.

Let

\begin{equation}
\label{eq34}
\overline y _j = \theta ^2y_j .
\end{equation}

Substituting (\ref{eq34}) into (\ref{eq33}), we get

\begin{equation}
\label{eq35}
\sum\limits_{j = 1}^m {a_{kj} \overline y _j = \theta ^2} .
\quad
k = 1,2,...,n
\end{equation}

Clearly, equations (\ref{eq35}) is a system of linear equations of $\overline y _j
$.

From (\ref{eq31}), we have

\begin{equation}
\label{eq36}
ml_j = \frac{r}{y_j }.
\end{equation}

From (\ref{eq17}), (\ref{eq20}), (\ref{eq34}) and (\ref{eq36}), we have

\begin{equation}
\label{eq37}
\sum\limits_{i = 1}^n {a_{ij} } c_i^2 = \frac{r\theta ^2}{m\overline y _j
}.
\end{equation}

And hence

\begin{equation}
\label{eq38}
\sum\limits_{i = 1}^n {a_{ij} (\frac{c_i }{\theta })^2} =
\frac{r}{m\overline y _j },
\quad
j = 1,2,...,m.
\end{equation}

Clearly, equations (\ref{eq38}) are also linear equations of $(c_i / \theta )^2$.
Equations (\ref{eq30}) are replaced by two system of linear equations, (\ref{eq35}) and
(\ref{eq38}). The problem seems to be very simple. Actually, there are many
combinations of (\ref{eq29}) and (\ref{eq30}),. Though, in these combinations, some of them
may not yield the solutions of this optimization problem, all the feasible
solutions of this problem are contained in the suitable combinations of
equations (\ref{eq29}) and equations (\ref{eq30}). If we solve the problem directly, then
there will be the same complexities as the simplex method for solving linear
programming problem. Furthermore, it is very difficulty to design the
algorithm and to implement the computer program if it is not impossible.
Even if the computer program is designed, then it might be a time-consuming
program. However, it can be concluded that the nonlinear optimization
problem is solvable theoretically. If the numbers $m$ and $n$ are very
small, say 3, then it is a simple problem to solve the systems of linear
equations. Generally, the extreme point of likelihood function is not
unique. The results of computer simulation show that the extreme point of
likelihood function seems to be unique. The computer simulations are
implemented when $m$ and $n$ are less than 10.

\textbf{Quantum theory approach}

In quantum mechanics, the wave function is linear combination of basis
functions and the normalization of the wave function requires that the sum
of the squares of coefficients should be unit. This gives us a clue to
remodel our problem and the problem becomes easier.

Let

\begin{equation}
\label{eq39}
\overline f (x) = \sum\limits_{i = 1}^n {u_i v_i } \varphi _i (x),
\end{equation}

\noindent
where $\varphi _i $ is the window function.

Let

\begin{equation}
\label{eq40}
\overline l = \prod\limits_{j = 1}^m {\overline f (x_j )}
\end{equation}

\noindent
be the likelihood function. The constraints are

\begin{equation}
\label{eq41}
\sum\limits_{i = 1}^n {u_i u_i } = r,
\end{equation}

\noindent
and

\begin{equation}
\label{eq42}
\sum\limits_{i = 1}^n {v_i v_i } = r.
\end{equation}

The problem is to maximize $\overline l $ subjected to constraints (\ref{eq41}) and
(\ref{eq42}). In order to find the connection of two models, we should define the
following notations.

Let

\begin{equation}
\label{eq43}
\sum\limits_{i = 1}^n {c_i } c_i = r.
\end{equation}

Let

\begin{equation}
\label{eq44}
l_j = \sum\limits_{i = 1}^n {c_i^2 \varphi _i (x_j )} .
\end{equation}

Let

\begin{equation}
\label{eq45}
l = \prod\limits_{j = 1}^m {l_j } .
\end{equation}

Let

\begin{equation}
\label{eq46}
S_n = \{(c_1 ,...,c_n );\sum\limits_{i = 1}^n {c_i } c_i = r\}.
\end{equation}

The first model is to find the extreme point of $l$ on $S_n $.

Let

\begin{equation}
\label{eq47}
\sum\limits_{i = 1}^n {u_i } u_i = r.
\end{equation}

Let

\begin{equation}
\label{eq48}
\sum\limits_{i = 1}^n {v_i } v_i = r.
\end{equation}

Let

\begin{equation}
\label{eq49}
\overline l _j = \sum\limits_{i = 1}^n {u_i v_i \varphi _i (x_j )} .
\end{equation}

Let

\begin{equation}
\label{eq50}
\overline l = \prod\limits_{j = 1}^m {\overline l _j } .
\end{equation}

Let

\begin{equation}
\label{eq51}
S_{2n} = \{(u_1 ,...,u_n ,v_1 ,...,v_n );\sum\limits_{i = 1}^n {u_i } u_i =
r,\sum\limits_{i = 1}^n {v_i } v_i = r\}.
\end{equation}

The second model is to find the extreme point of $\overline l $ on
$S_{2n} $. It is clear that $S_n $ and $S_{2n} $ are compact
subsets of $R^n$ and $R^{2n}$ respectively. Let $l$ and $\overline
l $ be the likelihood functions defined above. Clearly, both $l$
and $\overline l $ have maximum. Let $A$ be the set of all $l$.
Let $\overline A $ be the set of all $\overline l $. It is obvious
that $A \subseteq \overline A $. Therefore, the maximum of $A$ is
less than or equal to that of $\overline A $. It will be shown, in
theorem 1, that the extreme points of $\overline l $ should be
located at the points such that $u_i = v_i $, $i = 1,2,...,n$.
Therefore, the problem to maximize $l$ subjected to the constraint
(\ref{eq43}) is equivalent to that of maximizing $\overline l $
subjected to the constraints (\ref{eq47}) and (\ref{eq48}).

\textbf{Theorem 1. }For each observation $x_j $, if there
is\space$ \varphi _i $ such that $\varphi _i (x_j ) > 0$, then the
extreme points of $\overline l $ should be located at the points
such that $u_i = v_i $, $i = 1,2,...,n$.

\textbf{Proof. }We assume that

\begin{equation}
\label{eq52}
v_1 > u_1
\end{equation}

\noindent
and the maximum is $\overline l _M $, that is,

\begin{equation}
\label{eq53}
\overline l _M \ge \overline l
\end{equation}

\noindent
for all $\overline l $. Let

\begin{equation}
\label{eq54}
u'_i = v'_i = \theta \sqrt {u_i v_i } .
\end{equation}

By choosing a proper value $\theta $, constraints

\begin{equation}
\label{eq55}
\sum\limits_{i = 1}^n {u'_i } u'_i = r
\end{equation}

\noindent
and

\begin{equation}
\label{eq56}
\sum\limits_{i = 1}^n {v'_i } v'_i = r
\end{equation}

\noindent
are satisfied simultaneously. By Cauchy-Schwartz inequality and (\ref{eq52}), we get

\begin{equation}
\label{eq57}
\theta > 1
\end{equation}

\noindent
and hence we have

\begin{equation}
\label{eq58}
\overline l _M < \overline l
\end{equation}

\noindent
for some $\overline l $. This is a contradiction.

\textbf{The iteration procedures }

Tough we have stated and proved Theorem 1, we need a constructive procedure
to find the extreme point. It is not so easy to solve the nonlinear
optimization problem. Usually, the sequences are constructed by iteration
procedures. The well designed iteration procedures can generate monotonic
sequences which are useful in theory and application. With the nested
iteration procedures, the complicated problems such as mathematical
formulation, designing of the computation algorithm and the computer
programming can be solved in parallel. It seems that it is easier to
maximize $l$ than to maximize $\overline l $. The reason why we solve the
more complicated problem can be shown in the method of Lagrange's
multiplier. The strategy of solving the nonlinear optimization problem with
constraints is ignoring one of the constraints, say equation (\ref{eq48}). This can
be done by choosing the initial value of $v_i $, $v_i = \sqrt {r / n} $.
Then the optimization problem becomes simpler because only one constraint is
left. In order make it more clearly and precisely, we recall and define some
identities.

Let

\begin{equation}
\label{eq59}
\psi _i (x) = v_i \varphi _i (x).
\end{equation}

Let

\begin{equation}
\label{eq60}
\overline f (x) = \sum\limits_{i = 1}^n {u_i } \psi _i (x).
\end{equation}

Let

\begin{equation}
\label{eq61}
\overline l _j = \overline f (x_j ).
\end{equation}

Let

\begin{equation}
\label{eq62}
\overline l = \prod\limits_{j = 1}^m {\overline l _j } .
\end{equation}

The problem is to maximize $\overline l $ subjected to the constraint (\ref{eq47}).
After the values of $u_i $ being obtained, the value of $v_i $ is updated by
$\overline \theta \sqrt {u_i v_i } $. By choosing the factor $\overline
\theta $, the constraint (\ref{eq48}) is satisfied. Clearly, the iteration
procedures can be obtained. And the Cauchy-Schwartz inequality is able to
test the termination of the iteration procedures. First, we summarize the
whole procedures. Later, the associated mathematical theory of the procedure
will be shown. The procedures are:

\textbf{Step (i)}. \textit{Initialize the procedure by setting }$k = 1$\textit{ and }$v_i^k = \sqrt {r / n} , i = 1,2,...,n.$

\textbf{Step (ii)}\textit{, maximize }$\overline l $\textit{ subjected to the constraint (\ref{eq47}). Then values of }$u_i^k ,i = 1,2,...,n$\textit{, are obtained}.

\textbf{Step (iii)}, \textit{Check the condition }$\sum\limits_{i = }^n {u_i^k v_i^k + \varepsilon \ge r}
$\textit{ is satisfied or not, where }$\varepsilon $\textit{ is a small positive number to control the termination of the procedures.}

\textit{If the condition is satisfied, then stop the iteration procedures and the density estimator, }$\overline f (x) = \sum\limits_{i = 1}^n {u_i^k v_i^k } \varphi _i (x)$\textit{, is obtained. Otherwise, increase the value of }$k$\textit{ by one, set }$v_i^k
= \overline \theta ^k\sqrt {u_i^{k - 1} v_i^{k - 1} } $\textit{, here }$\overline \theta
^k$\textit{ is a factor to fit the constraint (\ref{eq48}). Then go to Step (ii) and proceed the procedures.}

\textbf{Remark 1}. \textit{From Cauchy-Schwartz inequality, the values of }$\overline \theta ^k$\textit{ must be greater than or equal to 1 and hence the set of the values of the likelihood function is an increasing sequence. }

Since step (i) and step (iii) are so simple, the only problem which is left
is how to complete the step (ii). Now, we will show how step (ii) can work
well. In order to complete the step(ii), another nested iteration procedures
will be designed and studied. In order to collaborate with the computer
algorithm, new notations must be introduced. Let $u_i^k $and $v_i^k $ be
obtained in the $k^{th}$ iteration. Let

\begin{equation}
\label{eq63}
\widehat{f}(x) = \sum\limits_{i = 1}^n {u_i^k (v_i^k } \varphi _i (x)).
\end{equation}

Let

\begin{equation}
\label{eq64}
\psi _i (x) = v_i^k \varphi _i (x).
\end{equation}

The simple notation,

\begin{equation}
\label{eq65}
\widehat{f}(x) = \sum\limits_{i = 1}^n {u_i } \psi _i (x)
\end{equation}

\noindent
shall be used hereafter.

\textbf{The constructive proof and the procedures of optimization}

\textbf{Lemma 1}. Let $\widehat{f}(x) = \sum\limits_{i = 1}^n {u_i
} \psi _i (x)$, where $\psi _i $ are nonnegative functions. Let
$\widehat{l}_j = \widehat{f}(x_j )$. Let $\widehat{l} = \prod\limits_{j =
1}^m {\widehat{l_j }} $ be the likelihood function. For each $x_j $, there
is a $\psi _i $ such that $\psi _i (x_j ) > 0$.Then there are constructive
procedures to maximize $\widehat{l}$ subjected to the constraint. We recall
the constraint (\ref{eq47})

\[
\sum\limits_{i = 1}^n {u_i u_i } = r.
\]

\textbf{Remark 2:} \textit{Since }$\psi _i $\textit{ are nonnegative functions and the constraint is invariant under the transformation, }$u_i = - u_i $\textit{, the solution of this optimization problem, }$u_i $\textit{, must be nonnegative.}

\textbf{Proof. }Let

\begin{equation}
\label{eq66}
b_{ij} = \psi _i (x_j ).
\end{equation}

Let

\begin{equation}
\label{eq67}
\widehat{l}_L = \widehat{l} - \lambda (\sum\limits_{i = }^n {u_i u_i } -
r).
\end{equation}

By the method of Lagrange's multiplier, we have

\begin{equation}
\label{eq68}
\frac{\partial \widehat{l_L }}{\partial u_k } = 0.
\end{equation}

By simple symbolic computation of derivatives, we get

\begin{equation}
\label{eq69}
\widehat{l}\sum\limits_{j = 1}^m {\frac{b_{kj} }{\widehat{l}_j }} - 2\lambda
u_k = 0,
\quad
k = 1,2,...,n.
\end{equation}

Multiplying (\ref{eq69}) by $u_k $ and taking the summation of index $k$, we get

\begin{equation}
\label{eq70}
\sum\limits_{k = 1}^n {\widehat{l}(\sum\limits_{j = 1}^m {\frac{b_{kj}
}{\widehat{l}_j }} - 2\lambda u_k )u_k } = 0,
\quad
k = 1,2,...,n.
\end{equation}

Interchanging the summations, we get

\begin{equation}
\label{eq71}
\sum\limits_{j = 1}^m {\widehat{l}(\sum\limits_{k = 1}^n {\frac{b_{kj} u_k
}{\widehat{l}_j }} - 2\lambda u_k u_k )} = 0.
\end{equation}

From (\ref{eq47}), (\ref{eq65}), (\ref{eq66}), (\ref{eq71}) and definition of $\widehat{l}_j $, we get

\begin{equation}
\label{eq72}
\widehat{l}m - 2\lambda r = 0,
\end{equation}

\noindent
and hence

\begin{equation}
\label{eq73}
\widehat{l} = (2\lambda r) / m.
\end{equation}

Substituting (\ref{eq73}) into (\ref{eq69}), we get

\begin{equation}
\label{eq74}
\sum\limits_{j = 1}^m {\frac{rb_{kj} }{m\widehat{l}_j }} - u_k = 0,
\quad
k = 1,2,...,n.
\end{equation}

Let ${\rm {\bf u}}$ and ${\rm {\bf b}}_j $be $n$ components vectors, where
${\rm {\bf u}} = \left[ {u_1 ,u_2 ,..,u_n } \right]^t$ and ${\rm {\bf b}}_j
= \left[ {b_{1j} ,b_{2j} ,...,b_{nj} } \right]^t$. By the constraint (\ref{eq47})
and the assumption of this lemma, ${\rm {\bf u}}$ and ${\rm {\bf b}}_j $ are
not zero vectors. Rewrite equations (\ref{eq74})

\begin{equation}
\label{eq75}
\frac{r}{m}\sum\limits_{j = 1}^m {\frac{{\rm {\bf b}}_j }{{\rm {\bf u}}
\cdot {\rm {\bf b}}_j }} - {\rm {\bf u}} = 0.
\end{equation}

Let

\begin{equation}
\label{eq76}
\alpha _k = \frac{1}{{\rm {\bf u}} \cdot {\rm {\bf b}}_k } \quad ,
\quad
k = 1,2,...,m.
\end{equation}

Substituting (\ref{eq76}) into (\ref{eq75}), we get

\begin{equation}
\label{eq77}
{\rm {\bf u}} = \frac{r}{m}\sum\limits_{j = 1}^m {\alpha _j {\rm {\bf b}}_j
} .
\end{equation}

Substituting (\ref{eq77}) into (\ref{eq76}), we obtain

\begin{equation}
\label{eq78}
\alpha _k = \frac{1}{\frac{r}{m}\sum\limits_{j = 1}^m {\alpha _j ({\rm {\bf
b}}_j \cdot {\rm {\bf b}}_k )} },
\quad
k = 1,2,...,m.
\end{equation}

Let

\begin{equation}
\label{eq79}
D_{ij} = \frac{r}{m}({\rm {\bf b}}_i \cdot {\rm {\bf b}}_j ),
\quad
i,j = 1,2,...,m.
\end{equation}

Substituting (\ref{eq79}) into (\ref{eq78}), we obtain

\begin{equation}
\label{eq80}
\sum\limits_{j = 1}^m {D_{kj} \alpha _k \alpha _j = 1} ,
\quad
k = 1,2,...,m.
\end{equation}

It should noticed that the major differences between (\ref{eq74}) and (\ref{eq80}) are the
range of the indices since $m$ and $n$ are different. If $m = 1$, then the
solution of equation (\ref{eq80}) can be obtained. From (\ref{eq77}), the lemma is proved.
Fortunately, if $m > 1$, then we can solve equations (\ref{eq80}) one by one. It is
very simple to show that the existence and the uniqueness of the solution of
(\ref{eq80}), we will complete the details of works in the following lemmas and
theorems. Now, we assume that the solution of (\ref{eq80}) can be obtained
effectively and the solution is unique. Therefore, the lemma is proved and
it seems that step (i), (ii) and (iii) can work well.

In deriving the equations, the systematic notations are adopted. Therefore,
the variables, $\alpha _k $ and $\alpha _j $ in the equations (\ref{eq80}) are
interchangeable. In order to simplify the problem, these equations will be
solved one by one in iteration procedures. The symmetry shall be destroyed
because only one variable, $\alpha _k $, will be focused. Usually, there are
at least two sets of variables in iteration procedures, one set is
associated with the old value and the other set is associated with the
updated new value. Therefore, we use the symbols with prime for new value.
In order to analyze the details of algorithm, the delta notation shall be
used, for example, $\alpha '_k = \alpha _k + \Delta \alpha _k $. Therefore,
there are different forms of equations (\ref{eq80}) in different notations. The
functions of different forms of equations (\ref{eq80}) are obvious because each form
is associated with a meaning. The error of each equation is denoted by $E_j
 \quad j = 1,2,...,m$, and $\Delta E_j $ is the variation of $E_j $ in iteration
procedure. The total sum of the absolute value of $E_j $ is denoted by $E$,
and $\Delta E$ is the variation of $E$. $E^{(j)}$ is the value of $E$ in
$j^{th}$ iteration.. These notations and their meanings shall be defined in
the context.

\textbf{Remark 3:} \textit{From identity (\ref{eq76}) and remark 2, }$\alpha _k $\textit{ must be nonnegative. Clearly, the solution of (\ref{eq77}), }${\rm {\bf u}}$\textit{, shall satisfy the constraint, }${\rm {\bf u}} \cdot {\rm {\bf
u}} = r$\textit{. It is not necessary to worry about that the quantity }${\rm {\bf u}} \cdot {\rm {\bf b}}_k $\textit{ in (\ref{eq76}) might be zero. The identity (\ref{eq76}) and (\ref{eq78}) are adopted for the convention of symbolic computations. These will be shown later.}

Though the likelihood function is highly nonlinear, equations in (\ref{eq80}) are a
system of quadratic equations. Intuitively, the solution of a single
quadratic equation can be obtained easily. In order to solve the equations
(\ref{eq80}) one by one, the nested iteration procedures are constructed. We write
one of them, say $k^{th}$ equation, the quadratic equation of $\alpha _k $,

\begin{equation}
\label{eq81}
D_{kk} \alpha _k^2 + (\sum\limits_{i \ne k}^n {D_{ik} \alpha _i )\alpha _k -
1 = 0} .
\end{equation}

Clearly, the only positive solution of (\ref{eq81}) is $( - s + \sqrt {s^2 + 4D_{kk}
} ) / (2D_{kk} ),$ where $s = \sum\limits_{i \ne k} {D_{ik} \alpha _i } $. If
the equations in (\ref{eq80}) can be solved one by one, then the problem becomes
simpler. Indeed, the equations in (\ref{eq80}) can be solved one by one and the sum
of all errors is reduced in each time. Basing on this fact, we are able to
design another set of iteration procedures step (a), (b) and (c) to solve
the problem. Now, we start to design the procedures.

Let

\begin{equation}
\label{eq82}
E_i = \sum\limits_{j = 1}^m {D_{ij} \alpha _i \alpha _j } - 1.
\end{equation}

Let

\begin{equation}
\label{eq83}
E = \sum\limits_{i = 1}^m {\left| {E_i } \right|} .
\end{equation}

In the iteration procedures, the values of $E_i $ and $E$ shall be changed.
Let $\Delta E_i $ be variation of $E_i $. Let $\Delta E$ be variation of
$E$. In order to collaborate with the algorithm, the nested iteration
procedures are designed in the step (ii). Clearly, the problem is to
minimize the value of $E$. And it must be proved that the minimum of $E$ is
zero. Therefore, the solution of (\ref{eq80}) and the solution of (\ref{eq77}) are obtained.
Now, we construct the iteration procedures to complete step(ii). The
associated mathematical lemmas and theorems of algorithm will emerge. First,
initialize the procedure by setting $\alpha _k = \sqrt {1 / (\overline D m}
)$, $k = 1,2,...,m$, where $\overline D $ is the maximum of $D_{ij} $. Then
the iteration procedures are:

\textbf{Step (a)}. \textit{Compute }$E_i = \sum\limits_{j = 1}^m {D_{ij} \alpha _i \alpha _j
} - 1,i = 1,2,...,m$\textit{, and }$E = \sum\limits_{i = 1}^m {\left| {E_i } \right|} $\textit{. Go to step (b).}

\textbf{Step (b)}. \textit{Test the condition whether }$E \le \delta $\textit{ is satisfied or not, where }$\delta $\textit{ is a small positive number to control the termination of the procedures. If }$E \le \delta $\textit{, then the desired results are obtained. Compute }$u_k $\textit{ by the identity (\ref{eq77}), }$k =
1,2,...,n$\textit{, and terminate the iteration. Otherwise, go to step (c).}

\textbf{Step (c)}, \textit{Find the largest element of the set of all }$\left| {E_i } \right|$\textit{. Suppose that the largest element is }$\left| {E_k } \right|$\textit{ for some }$k$\textit{. Eliminate }$E_k $\textit{ by updating the value of }$\alpha
_k $\textit{ by }$\alpha '_k , \alpha '_k = ( - s + \sqrt {s^2 + 4D_{kk} } ) / (2D_{kk}
)$\textit{, where }$s = \sum\limits_{i \ne k} {D_{ik} \alpha _i } $\textit{. Go to Step (a).}

Now, there will be no difficulty to implement steps (a), (b) and (c).
Intuitively, steps (a), (b) and (c) shall be terminated in finite steps if
the values of $E$ is strictly decreasing sequence which converges to zero.
In lemma 3, it will be proved that the values of $E$ is a decreasing
sequence. Lemma 2 will support lemma 3$.$ In lemma 5, it will be proved that
the values of $E$ is a strictly decreasing sequence which converges to zero.
Lemma 4 will support lemma 5.

\textbf{Lemma 2}. All iterations, steps (a), (b) and (c), the set of all
$\alpha _k $, $k = 1,2...,m$, are bounded above and the set of all $\alpha
_k $, $k = 1,2...,m$, are bounded below by a positive number , say $B$, $B >
0$. That is, $\alpha _k > B$, $k = 1,2...,m$.

\textbf{Remark 4}. \textit{What we mean all }$\alpha _k $\textit{ is including all }$\alpha _k $\textit{ and all }$\alpha '_k .$

\textbf{Proof.} The value of $\alpha _k $ is either the initial value $\sqrt
{1 / (\overline D m} )$ or the updated value $( - s + \sqrt {s^2 + 4D_{kk} }
) / (2D_{kk} )$. It is very easy to verify the following inequalities

\begin{equation}
\label{eq84}
\frac{ - s + \sqrt {s^2 + 4D_{kk} } }{2D_{kk} } < \frac{ - s + \sqrt {s^2 +
2sD_{kk} + D_{kk}^2 } }{2D_{kk} } = \frac{1}{2},
\end{equation}

\noindent
when $s \ge 2$.

\begin{equation}
\label{eq85}
\frac{ - s + \sqrt {s^2 + 4D_{kk} } }{2D_{kk} } < \frac{\sqrt {4 + 4D_{kk} }
}{2D_{kk} },
\end{equation}

\noindent
when $s < 2$.

It is obvious that $\alpha _k $ are bounded above. Since $s = \sum\limits_{i
\ne k} {D_{ik} \alpha _i } $, s is bounded above. Next, we are going to
prove that there is a positive number $B$ such that $\alpha _k > B$, $k =
1,2...,m$, in all iterations. Clearly, $\alpha _k $ are either the initial
value or updated by $( - s + \sqrt {s^2 + 4D_{kk} } ) / (2D_{kk} )$. The
derivative of $( - s + \sqrt {s^2 + 4D_{kk} } ) / (2D_{kk} )$ is $( - 1 +
s / \sqrt {s^2 + 4D_{kk} } ) / (2D_{kk} )$, which is negative for all $s \ge
0$. Therefore, $( - s + \sqrt {s^2 + 4D_{kk} } ) / (2D_{kk} )$ is a decreasing
function of $s$. It is obvious that

\begin{equation}
\label{eq86}
\mathop {\lim }\limits_{s \to \infty } ( - s + \sqrt {s^2 + 4D_{kk} } ) /
(2D_{kk} ) = 0.
\end{equation}

Since $s$ is bounded above, $( - s + \sqrt {s^2 + 4D_{kk} } ) / (2D_{kk} )$
has a positive lower bound,. Therefore, $\alpha _k $ is bounded below by a
positive lower bound, say$B$.

\textbf{Remark 5.} \textit{Lemma 2 does not imply }$u_k $\textit{ are bounded below by a positive number, some }$u_k $\textit{ might tend to zero.}

\textbf{Lemma 3. }The values of $E$ in iteration procedures, step (a), (b)
and (c), is a decreasing sequence.

\textbf{Proof.} From (\ref{eq82}) and (\ref{eq83}), we find that equations (\ref{eq80}) can be
solved one by one. One of equations (\ref{eq80}) with one variable, say $\alpha '_k
$, will be solved. The error of the equation with index $k$, $E_k $, is
removed completely in step (c). Therefore,

\begin{equation}
\label{eq87}
\left| {\Delta E_k } \right| = \left| {E_k } \right|
\end{equation}

\noindent
for the particular index $k$ and it might be that

\begin{equation}
\label{eq88}
\left| {\Delta E_j } \right| \ne \left| {E_j } \right|
\end{equation}

\noindent
when $j \ne k$. Though there are two roots of a quadratic equation, only one
of them is positive. From equation (\ref{eq81}), it must be $( - s + \sqrt {s^2 +
4D_{kk} } ) / (2D_{kk} )$.

Let

\begin{equation}
\label{eq89}
\Delta \alpha _k = \alpha '_k - \alpha _k .
\end{equation}

The value of $\Delta \alpha _k $ is the difference of two positive numbers
which are bounded above. Clearly,

\begin{equation}
\label{eq90}
\left| {\Delta \alpha _k } \right| \le \alpha _k
\end{equation}

\noindent
when $\Delta \alpha _k \le 0$. In the step (c), the value of $E$ is reduced
by $\Delta E$. In order to update the value of $\alpha _k $, we rewrite the
equation (\ref{eq81})

\begin{equation}
\label{eq91}
\alpha '_k \sum\limits_{i \ne k}^m {D_{ki} \alpha _i } + D_{kk} \alpha '_k
\alpha '_k - 1 = 0.
\end{equation}

Some times, it is more convenient to use the delta notation. Therefore,
equation (\ref{eq91}) becomes

\begin{equation}
\label{eq92}
(\alpha _k + \Delta \alpha _k )\sum\limits_{i \ne k}^m {D_{ki} \alpha _i } +
D_{kk} (\alpha _k + \Delta \alpha _k )^2 - 1 = 0,
\end{equation}

\begin{equation}
\label{eq93}
\Delta \alpha _k \sum\limits_{i = 1}^m {D_{ki} \alpha _i } + \Delta \alpha
_k D_{kk} \alpha _k + D_{kk} (\Delta \alpha _k )^2 + \alpha _k
\sum\limits_{i = 1}^m {D_{ki} \alpha _i } - 1 = 0.
\end{equation}

From (\ref{eq82}), we get

\begin{equation}
\label{eq94}
\Delta E_k = \alpha _k \sum\limits_{i = 1}^m {D_{ki} \alpha _i } - 1,
\end{equation}

\noindent
for this particular index $k$.

Rewrite (\ref{eq93})

\begin{equation}
\label{eq95}
\Delta \alpha _k \sum\limits_{i = 1}^m {D_{ki} \alpha _i } + \Delta \alpha
_k D_{kk} \alpha _k + D_{kk} (\Delta \alpha _k )^2 = - (\alpha _k
\sum\limits_{i = 1}^m {D_{ki} \alpha _i } - 1).
\end{equation}

Clearly,

\begin{equation}
\label{eq96}
\left| {\Delta \alpha _k \sum\limits_{i = 1}^m {D_{ki} \alpha _i } + \Delta
\alpha _k D_{kk} \alpha _k + D_{kk} (\Delta \alpha _k )^2} \right| = \left|
{\Delta E_k } \right|,
\end{equation}

\begin{equation}
\label{eq97}
\left| {\Delta \alpha _k } \right|\left| {(\sum\limits_{i \ne k}^m {D_{ki}
\alpha _i } + D_{kk} 2\alpha _k + D_{kk} \Delta \alpha _k )} \right| =
\left| {\Delta E_k } \right|.
\end{equation}

All quantities in $(\sum\limits_{i \ne k}^m {D_{ki} \alpha _i } + D_{kk}
2\alpha _k + D_{kk} \Delta \alpha _k )$, except $\Delta \alpha _k $, are
positive. From (\ref{eq90}), for any case,

\begin{equation}
\label{eq98}
\left| {(\sum\limits_{i \ne k}^m {D_{ki} \alpha _i } + D_{kk} \alpha _k )}
\right| \le \left| {(\sum\limits_{i \ne k}^m {D_{ki} \alpha _i } + D_{kk}
2\alpha _k + D_{kk} \Delta \alpha _k )} \right|.
\end{equation}

Therefore,

\begin{equation}
\label{eq99}
\left| {\Delta \alpha _k } \right|\left| {(\sum\limits_{i \ne k}^m {D_{ki}
\alpha _i } + D_{kk} \alpha _k )} \right| \le \left| {\Delta E_k } \right|.
\end{equation}

If we write whole system of equations (\ref{eq80}), then the upper bound of all
$\left| {\Delta E_i } \right|$, $i \ne k$, can be figured out. From (\ref{eq82}), we
get

\begin{equation}
\label{eq100}
\Delta E_i = (\alpha _i + \Delta \alpha _i )\sum\limits_{j = 1}^m {D_{ij}
(\alpha _j + \Delta \alpha _j )} - 1 - E_i .
\end{equation}

\noindent
for all $i$. But

\begin{equation}
\label{eq101}
\Delta \alpha _i = 0,
\quad
i \ne k.
\end{equation}

From (\ref{eq100}) and (\ref{eq101}), we get

\begin{equation}
\label{eq102}
\Delta E_i = \Delta \alpha _k D_{ik} \alpha _i ,
\quad
i \ne k.
\end{equation}

Since

\begin{equation}
\label{eq103}
D_{ik} = D_{ki} ,
\end{equation}

\begin{equation}
\label{eq104}
\Delta \alpha _k D_{ik} \alpha _i = \Delta \alpha _k D_{ki} \alpha _i .
\end{equation}

From the inequality (\ref{eq99}) and (\ref{eq102}), we get

\begin{equation}
\label{eq105}
\sum\limits_{i \ne k} {\left| {\Delta E_i } \right| + \left| {\Delta \alpha
_k } \right|D_{kk} \alpha _k < \left| {\Delta E_k } \right|} \quad .
\end{equation}

From (\ref{eq105}),

\begin{equation}
\label{eq106}
\left| {\Delta \alpha _k } \right|D_{kk} \alpha _k < \left| {\Delta E_k }
\right| - \sum\limits_{i \ne k} {\left| {\Delta E_i } \right|} .
\end{equation}

From (\ref{eq83}), we get

\begin{equation}
\label{eq107}
\Delta E = \sum\limits_{i = 1}^m {\Delta \left| {E_i } \right|} .
\end{equation}

Since

\begin{equation}
\label{eq108}
\left| {a + b} \right| \ge \left| a \right| - \left| b \right|
\end{equation}

\noindent
for any $a$ and $b$,

\begin{equation}
\label{eq109}
\left| {\Delta E} \right| \ge \left| {\Delta E_k } \right| - \sum\limits_{i
\ne k} {\left| {\Delta E_i } \right|} .
\end{equation}

Clearly, $\Delta E$ is negative and dominated by $\left| {\Delta E_k }
\right|$,

\begin{equation}
\label{eq110}
\left| {\Delta E} \right| > \left| {\Delta \alpha _k } \right|D_{kk} \alpha
_k .
\end{equation}

An hence the set of the values of $E$ generated by iterations is a
decreasing sequence, We have proved the lemma.

For each iteration, the value of $E$ is denoted by a symbol, say $E^{(i)}$
in the $i^{th}$ iteration. The notations $E_i $ and $E^{(i)}$ are associated
with different meanings. Let $\lim _{i \to \infty } E^{(i)} = E^\infty $.
Clearly, the lower bound of $\left| {\Delta \alpha _k } \right|D_{kk} \alpha
_K $ will serve for two purposes, one is to prove that the sequence
$E^{(i)}$ is a strictly decreasing sequence and the other is to prove that
$E^\infty = 0$.

\textbf{Lemma 4.} If $E^\infty > 0$ and $k$ is the index such that $\left|
{E_k } \right| \ge \left| {E_i } \right| \quad i = 1,2,...,m$, then the set of
all $\left| {\Delta \alpha _k } \right|D_{kk} \alpha _k $, in all iterations
of step (a), (b) and (c) has a nonzero lower bound.

\textbf{Proof.} It is obvious that

\begin{equation}
\label{eq111}
\left| {E_k } \right| \ge E^\infty / m.
\end{equation}

In each iteration procedure, only one equation is solved. From (\ref{eq97}) and
(\ref{eq111}), we get

\begin{equation}
\label{eq112}
\left| {\Delta \alpha _k } \right|\left| {(\sum\limits_{i \ne k}^m {D_{ki}
\alpha _i } + D_{kk} 2\alpha _k + D_{kk} \Delta \alpha _k )} \right| \ge
\frac{E^\infty }{m}.
\end{equation}

The first term absorbing $D_{kk} \alpha _k $ from the second term , we get

\begin{equation}
\label{eq113}
\left| {\Delta \alpha _k } \right|\left| {\sum\limits_{i = 1}^m {D_{ki}
\alpha _i } + D_{kk} \alpha _k + D_{kk} \Delta \alpha _k } \right| \ge
\frac{E^\infty }{m}.
\end{equation}

Since $\left| {\Delta \alpha _k } \right|$ and $\alpha _i $ are bounded
above, $\left| {\sum\limits_{i = 1}^m {D_{ki} \alpha _i } + D_{kk} \alpha _k
+ D_{kk} \Delta \alpha _k } \right|$ is also bounded above, say

\begin{equation}
\label{eq114}
\left| {\sum\limits_{i = 1}^m {D_{ki} \alpha _i } + D_{kk} \alpha _k +
D_{kk} \Delta \alpha _k } \right| < M.
\end{equation}

From (\ref{eq113}) and (\ref{eq114}), we get

\begin{equation}
\label{eq115}
\left| {\Delta \alpha _k } \right| > E^\infty / (Mm).
\end{equation}

Therefore, $\left| {\Delta \alpha _k } \right|$is bounded below by a
positive number and hence $\left| {\Delta \alpha _k } \right|D_{kk} \alpha
_k $ is bounded below by a positive number in all iterations. Therefore, we
have proved the lemma.

\textbf{Lemma 5. }$\mathop {\lim }\limits_{k \to \infty } E^k = 0$, that is,
$E^\infty = 0$.

\textbf{Proof}. For any $\varepsilon $, $\varepsilon > 0$, there is an
positive integer $N$ such that

\begin{equation}
\label{eq116}
E^{(j)} - \varepsilon < E^\infty ,
\end{equation}

\noindent
whenever $j \ge N$. Since $\left| {E_k } \right|$ is the largest one in the
$j^{th}$iteration,

\begin{equation}
\label{eq117}
\left| {E_k } \right| \ge E^\infty / m.
\end{equation}

From inequality (\ref{eq110}),

\begin{equation}
\label{eq118}
E^{(j + 1)} + \left| {\Delta \alpha _k } \right|D_{kk} \alpha _k < E^{(j)}.
\end{equation}

Therefore,

\begin{equation}
\label{eq119}
E^{(j + 1)} + \left| {\Delta \alpha _k } \right|D_{kk} \alpha _k -
\varepsilon < E^{(j)} - \varepsilon .
\end{equation}

If we assume that

\begin{equation}
\label{eq120}
E^\infty > 0.
\end{equation}

By lemma 4, the set of all $\left| {\Delta \alpha _k } \right|D_{kk} \alpha
_k $ has a nonzero lower bound. We choose $\varepsilon $ such that
$\varepsilon $ is less than the lower bound of $\left| {\Delta \alpha _k }
\right|D_{kk} \alpha _k $. That is,

\begin{equation}
\label{eq121}
\left| {\Delta \alpha _k } \right|D_{kk} \alpha _k - \varepsilon > 0.
\end{equation}

Then

\begin{equation}
\label{eq122}
E^{(j + 1)} < E^{(j)} - \varepsilon .
\end{equation}

From (\ref{eq116}), we get

\begin{equation}
\label{eq123}
E^{(j + 1)} < E^\infty .
\end{equation}

It is a contradiction because $E^{(j)} \ge E^\infty $ for all $j$.
Therefore, we have proved the lemma and hence $E^\infty = 0$.

Since $E^\infty = 0$, the iteration procedures, step (a), step (b) and step
(c), should terminate in finite steps of iterations and step (ii) can be
executed completely. Therefore, lemma 1 is proved completely. In lemma 7, it
will be proved that the iteration procedures, step (i), step(ii) and step
(iii), shall be terminated in finite steps. Lemma 6 will support lemma 7$.$

\textbf{Lemma 6}. Let $\overline \theta ^k$ be obtained in the iteration
procedures, step (i), step(ii) and step (iii). Then $\mathop {\lim
}\limits_{k \to \infty } \overline \theta ^k = 1$.

\textbf{Proof}. It is obvious that

\begin{equation}
\label{eq124}
\overline \theta ^k \ge 1.
\end{equation}

\noindent
and hence $\widehat{l}^k$ is an increasing sequence.

Let

\begin{equation}
\label{eq125}
\mathop {\lim }\limits_{k \to \infty } \widehat{l}^k = \widehat{l}^\infty .
\end{equation}

For any $e > 0$, there is $\widehat{l}^k$ such that

\begin{equation}
\label{eq126}
\widehat{l}^k + e > \widehat{l}^\infty .
\end{equation}

If $\mathop {\lim }\limits_{k \to \infty } \overline \theta ^k$does not
exist, then there exist $\varepsilon > 0$, for any $K$, there is $k > K$
such that

\begin{equation}
\label{eq127}
\overline \theta ^k > 1 + \varepsilon .
\end{equation}

From (\ref{eq64}), (\ref{eq65}) and the definition $\widehat{l}$, we get

\begin{equation}
\label{eq128}
\widehat{l}^{k + 1} > (\overline \theta ^k)^m\widehat{l}^k,
\end{equation}

\noindent
where $m$ is the sample size. Since

\begin{equation}
\label{eq129}
(1 + \varepsilon )^m > 1 + m\varepsilon ,
\end{equation}

\begin{equation}
\label{eq130}
\widehat{l}^{k + 1} > (1 + m\varepsilon )\widehat{l}^k.
\end{equation}

Therefore,

\begin{equation}
\label{eq131}
\widehat{l}^{k + 1} > \widehat{l}^k + m\varepsilon \widehat{l}^1.
\end{equation}

Choosing $e = m\varepsilon \widehat{l}^1$, we have

\begin{equation}
\label{eq132}
\widehat{l}^{k + 1} > \widehat{l}^k + e.
\end{equation}

From (\ref{eq126}), we get

\begin{equation}
\label{eq133}
\widehat{l}^{k + 1} > \widehat{l}^\infty .
\end{equation}

It is a contradiction. Therefore,

\begin{equation}
\label{eq134}
\mathop {\lim }\limits_{k \to \infty } \overline \theta ^k = 1.
\end{equation}

\textbf{Lemma 7}. Let $P^k = \sum\limits_{i = 1}^n {u_i^k } v_i^k $. Then
$\mathop {\lim }\limits_{k \to \infty } P^k = r$ and the iteration
procedures, step (i), step(ii) and step (iii), shall be terminated in finite
steps.

\textbf{Proof.} From the definition of $\overline \theta ^k$in step (iii),
we get

\begin{equation}
\label{eq135}
\overline \theta ^k\overline \theta ^k\sum\limits_{i = 1}^n {u_i^{k - 1} }
v_i^{k - 1} = r.
\end{equation}

From (\ref{eq134}) and (\ref{eq135}), we get

\begin{equation}
\label{eq136}
\mathop {\lim }\limits_{k \to \infty } P^k = r.
\end{equation}

Therefore, the iteration procedures, step (i), step(ii) and step (iii),
shall be terminated in finite steps.

Lemma 8 will show the result of theorem 1 can be obtained by constructive
method.

\textbf{Lemma 8.} Let $w_i^k = \left| {u_i^k - v_i^k } \right|$, $i =
1,2,...,n$, $k = 1,2,...,$ be a set sequences generated by the iteration
procedures, step (i), step(ii) and step (iii). Then $\mathop {\lim
}\limits_{k \to \infty } w_i^k = 0$, $i = 1,2,...,n$.

\textbf{Remark 6.}\textit{ By Cauchy-Schwartz inequality, }$(\sum\limits_{i = 1}^n {u_i^k v_i^k )^2 \le }
\sum\limits_{i = 1}^n {(u_i^k )^2} \sum\limits_{i = 1}^n {(v_i^k )^2} $\textit{, the equal sign hold only if }$u_i^k
= v_i^k , i = 1,2,...,n$\textit{. Intuitively, it is obvious that the condition in step (iii) must be satisfied. Otherwise, }$\overline l $\textit{ is not bounded above and }$\overline l $\textit{ does not have maximum. }

\textbf{Proof.} By simple computation, we get

\begin{equation}
\label{eq137}
\sum\limits_{i = 1}^n {(u_i^k - v_i^k )^2 = \sum\limits_{i = 1}^n {u_i^k } }
u_i^k + \sum\limits_{i = 1}^n {v_i^k } v_i^k - 2\sum\limits_{i = 1}^n {u_i^k
} v_i^k .
\end{equation}

In all iteration procedures, the constraints (\ref{eq47}) and (\ref{eq48}) must be
satisfied. Therefore

\begin{equation}
\label{eq138}
\sum\limits_{i = 1}^n {(u_i^k - v_i^k )^2 = } 2r - 2\sum\limits_{i = 1}^n
{u_i^k } v_i^k .
\end{equation}

From (\ref{eq136}) and (\ref{eq138}), we get

\begin{equation}
\label{eq139}
\mathop {\lim }\limits_{k \to \infty } w_i^k = 0,
\quad
i = 1,2,...,n.
\end{equation}

Combining the results theorem 1 and lemma 8, the problem of optimization is
solved almost.

\textbf{The unique theorem}

\textbf{Theorem 2}. The solution of equations (\ref{eq80}) is unique.

\textbf{Proof}: Of course, only the positive solutions make sense. Let ${\rm
{\bf e}}_i = \alpha _i {\rm {\bf b}}_i $, where $\alpha _i $ is a solution
that we have obtained by the iteration procedures step(a), step(b) and step
(c) . Let

\begin{equation}
\label{eq140}
e_{ij} = {\rm {\bf e}}_i \cdot {\rm {\bf e}}_j .
\end{equation}

From (\ref{eq80}) and (\ref{eq140}), we get

\begin{equation}
\label{eq141}
\sum\limits_{j = 1}^m {e_{ij} = 1} ,
\quad
i = 1,2...,m.
\end{equation}

Consider the following equations,

\begin{equation}
\label{eq142}
\sum\limits_{j = 1}^m {e_{ij} \beta _i \beta _j = 1} ,
\quad
i = 1,2...,m.
\end{equation}

Here $\beta _i $, $i = 1,2...,m$, are unknowns. Then

\begin{equation}
\label{eq143}
\beta _i = 1,
\quad
i = 1,2...,m,
\end{equation}

\noindent
is a solution of (\ref{eq141}). If the there is another solution set, say

\begin{equation}
\label{eq144}
\beta _1 \ge \beta _2 \ge ,... \ge \beta _m .
\end{equation}

From (\ref{eq141}), we know that the equal sign can not hold all times. From (\ref{eq142}),
we have

\begin{equation}
\label{eq145}
\beta _1 \sum\limits_{j = 1}^m {e_{1j} \beta _j = 1} ,
\end{equation}

\noindent
and

\begin{equation}
\label{eq146}
\beta _m \sum\limits_{j = 1}^m {e_{mj} \beta _j = 1} .
\end{equation}

It is obvious that

\begin{equation}
\label{eq147}
e_{ii} > 0,
\end{equation}

\noindent
for all $i$. Therefore,

\begin{equation}
\label{eq148}
\beta _1 \sum\limits_{j = 1}^m {e_{1j} \beta _j > } \beta _1 \sum\limits_{j
= 1}^m {e_{1j} \beta _m = } \beta _1 \beta _m \sum\limits_{j = 1}^m {e_{1j}
} > \beta _m \sum\limits_{j = 1}^m {e_{mj} } \beta _j .
\end{equation}

We have used the identities (\ref{eq141}) at least two times. From (\ref{eq142}), (\ref{eq146}) and
(\ref{eq148}), we find that it is a contradiction. We have completed the proof the
theorem.

\textbf{Theorem 3}. The solution of (\ref{eq74}) is unique and hence the maximum
value obtained by step (ii) is the global maximum on the sphere
$\sum\limits_{i = 1}^n {u_i u_i } = r$.

Proof. For simplicity, we use (\ref{eq75}), the vector notations, instead of (\ref{eq74}).
If there are two solutions say ${\rm {\bf u}}$ and ${\rm {\bf u}}'$.
Therefore,

\begin{equation}
\label{eq149}
\frac{r}{m}\sum\limits_{j = 1}^m {\frac{{\rm {\bf b}}_j }{{\rm {\bf u}}
\cdot {\rm {\bf b}}_j }} - {\rm {\bf u}} = 0,
\end{equation}

And

\begin{equation}
\label{eq150}
\frac{r}{m}\sum\limits_{j = 1}^m {\frac{{\rm {\bf b}}_j }{{\rm {\bf u}}'
\cdot {\rm {\bf b}}_j }} - {\rm {\bf u}}' = 0.
\end{equation}

Let

\begin{equation}
\label{eq151}
\alpha _k = \frac{1}{{\rm {\bf u}} \cdot {\rm {\bf b}}_k },
\quad
k = 1,2,...,m.
\end{equation}

Let

\begin{equation}
\label{eq152}
\alpha '_k = \frac{1}{{\rm {\bf u}}' \cdot {\rm {\bf b}}_k },
k = 1,2,...,m.
\end{equation}

Substituting (\ref{eq151}) into (\ref{eq149}), we get

\begin{equation}
\label{eq153}
{\rm {\bf u}} = \frac{r}{m}\sum\limits_{j = 1}^m {\alpha _j {\rm {\bf b}}_j
} ,
\end{equation}

Substituting (\ref{eq152}) into (\ref{eq150}), we get

\begin{equation}
\label{eq154}
{\rm {\bf u}}' = \frac{r}{m}\sum\limits_{j = 1}^m {\alpha '_j {\rm {\bf
b}}_j } .
\end{equation}

Substituting (\ref{eq153}) into (\ref{eq151}), we get

\begin{equation}
\label{eq155}
\alpha _k = \frac{1}{\frac{r}{m}\sum\limits_{j = 1}^m {\alpha _j ({\rm {\bf
b}}_j \cdot {\rm {\bf b}}_k )} },
\quad
k = 1,2,...,m,
\end{equation}

Substituting (\ref{eq154}) into (\ref{eq152}), we get

\begin{equation}
\label{eq156}
\alpha '_k = \frac{1}{\frac{r}{m}\sum\limits_{j = 1}^m {\alpha '_k ({\rm
{\bf b}}_j \cdot {\rm {\bf b}}_k )} },
\quad
k = 1,2,...,m.
\end{equation}

From (\ref{eq79}), (\ref{eq155}) and (\ref{eq156}), we get two systems of equations

\begin{equation}
\label{eq157}
\sum\limits_{j = 1}^m {D_{kj} \alpha _k \alpha _j = 1} ,
\quad
k = 1,2,...,m,
\end{equation}

\noindent
and

\begin{equation}
\label{eq158}
\sum\limits_{j = 1}^m {D_{kj} \alpha '_k \alpha '_j = 1} ,
\quad
k = 1,2,...,m.
\end{equation}

By theorem 2.

\begin{equation}
\label{eq159}
\alpha '_k = \alpha _k ,
\end{equation}

\noindent
for $k = 1,2,...,m$. From (\ref{eq153}), (\ref{eq154}) and (\ref{eq159}), we get

\begin{equation}
\label{eq160}
{\rm {\bf u}}' = {\rm {\bf u}}.
\end{equation}

And hence the maximum which is obtained in this algorithm is the global
maximum on the manifold, the sphere $\sum\limits_i^n {u_i } u_i = r$.

\textbf{Numerical examples}

No matter how good might the paper be, the final result must be verified by
numerical examples. There are three examples. The results are shown in
Figure 1, Figure 2 and Figure 3. The estimator is obtained by Bernstein
polynomials [2], [3].

Let

\begin{equation}
\label{eq161}
\varphi _i (x) = N_i (n! / (i!(n - i)!))x^i(1 - x)^{n - i},
\quad
0 \le x \le 1.
\end{equation}

Here, $N_i $ is a factor to make

\begin{equation}
\label{eq162}
\int_0^1 {\varphi _i (x)dx = 1} .
\end{equation}

Let

\begin{equation}
\label{eq163}
\widehat{f}(x) = \sum\limits_{i = 1}^n {c_i \varphi _i (x)} .
\end{equation}

We use the density estimator which is defined in the very beginning identity
(\ref{eq6}) though it is computed by (\ref{eq39}). All the observations, $x_1 $, $x_2
$,\ldots , $x_m $, must be contained in an interval [a, b]. It is a simple
work to transform the interval [a, b] to the interval [0, 1].

Example 1.

The density function is defined on $[0,\infty )$,

\[
f(x) = \exp ( - x).
\]

Example 2.

The density function is defined on $[0,4]$,

\[
f(x) = 2 / 3,
\]

\noindent
when $1 \le x \le 2$;

\[
f(x) = 1 / 3,
\]

\noindent
when $3 \le x \le 4$;

\[
f(x) = 0,
\]

\noindent
otherwise.

Example 3.

The density function is defined on $[0,4]$,

\[
f(x) = 1,
\]

\noindent
when $0 \le x \le 1 / 2$;

\[
f(x) = 1 / 2,
\]

\noindent
when $1 \le x \le 3 / 2$;

\[
f(x) = 1 / 2,
\]

\noindent
when $3 \le x \le 7 / 2$;

\[
f(x) = 0,
\]

\noindent
otherwise.

The density function of Example 2 and Example 3 are not continuous and hence
it is inappropriate to apply Bernstein polynomial to these examples. If the
piecewise spline is used, then the result shall be better actually. We will
not discuss the piecewise Bernstein polynomial in this paper. Comparing with
the existent method [11],[12] etc., the spline kernel or spline window is a
new method with potential because there will be new useful splines that
might be designed in near future. At least, there are three useful splines,
B-spline, Cubic spline and Bezier spline. The works of source program
designing, debugging and maintaining are more difficult than the
mathematical proofs because they are tedious works. More than four kernel
functions or window functions are tested, including B-spline, overlap
B-spline, Bezier spline and piecewise Bezier spline. Though we do not show
the result of B-spline approach, most programs are tested by B-spline method
first. We will not list the definition of B-spline because it is available
to find the definition of the splines in the books of numerical analysis. It
seems that B-spline method can be taken as the priori in Bayesian approach
and hence piecewise Bezier spline method can be taken as posterior in
Bayesian approach. Unlike the Bezier spline, the B-spline need the extra
control points, the knot points [2], and these knot points make the programs
more complicated and difficult. In the testing program, there are about 300
window functions are used in B-spline method. It is a good experiment to
solve about 300 nonlinear equations. The whole work is accomplished by using
the oldest fashion and the most modern language, visual fortran. If it is
necessary, then the fortran source programs will be appended.

\textbf{Discussion and conclusion. }

The algorithm is so attractive that it is not necessary to prove that these
sequences $u_i^k $, $v_i^k $ and $u_i^k v_i^k $, $i = 1,2,...,n$ converge.
The results of computer output show that these sequences $u_i^k $, $v_i^k $,
$u_i^k v_i^k $ converge . Moreover, $\sum\limits_{i = 1}^n {u_i^k v_i^k } $
is an increasing sequence and $\overline \theta ^k$ is a decreasing
sequence. The algorithm is to maximize likelihood function $l$ and to
terminate the procedures by the condition $\sum\limits_{i = 1}^n {u_i v_i }
> r - \varepsilon $. The constraints, (\ref{eq47}) and (\ref{eq48}), are satisfied in every
step. Since it has been proved that $\mathop {\lim w_i^k =
0}\limits_{k \to \infty } $, all $\left| {u_i - v_i } \right|$ are
very small when the iteration procedures are terminated. In this
paper, we do not prove the convergence of the sequences $u_i^k $,
$v_i^k $ and $u_i^k v_i^k $, $i = 1,2,...,n$. It should be
reminded that the problem is to maximize the likelihood function
subjected to the conditions $\sum\limits_{i = 1}^n {u_i } u_i = r$
and $v_i = u_i $, $i = 1,2,...,n$. We think that the problem is
solved almost. It is still an open problem whether the iterations
procedures, step (i), (ii) and (iii), will serve the purpose or
not, for finding the global maximum of $\overline A $? Of course,
step (i) play important role for searching for the global maximum
of $\overline A $, it seems to be so. We think that only if the
initial value of $v_i $ in step (i) is set $v_i \ne 0$ for all
$i$, then the procedures will find the global maximum. But the
proof is not completed yet. It is the unique theorems, theorem 2
and theorem 3, that simplify the complicated problem and gives us
the motivation to prove the global property. Since the consistency
of parametric estimator has been proved statistician [1], the only
problem left is finding the point which will yield the global
maximum of likelihood function. To the best knowledge of the
authors, there is no definite answer for finding the global
maximum of nonlinear optimization problems. Though the problem do
not be solved completely in theory, the work and its related
algorithm are very useful in practical problem.

The proof of the lemma 8 is short and simple because this is the final
version. The first version is abandoned because it is lengthy and
complicated. In the first version of the proof, we use the method of
variation. The technique of the first proof in lemma 8 is almost the same as
that of quantum physics, especially in quantum field theory and string
theory [13].

To follow the approach of most nonparametric approaches, we use the
advantage the probability amplitude which is introduced in the quantum
theory. Though the orthogonal polynomials are also used in both
nonparametric approaches and quantum theory, we use the Bernstein
polynomial. It is the Bersnstein polynomials that unify and simplify
fundamental problems such as parametric approach and nonparametric approach,
consistency of the estimator and the most difficult problem of density
estimation, the nonlinear optimization problem. The probability amplitude is
stressed most books of quantum physics [14].

It should be clarified that the research work is initiated and completed
finally by Yeong-Shyeong Tsai. Without the consultation with Lu-Hsing Tsai,
Hung-Ming Tsai and Po-Yu Tsai in quantum physics and personal computing
system, and the consultation with Yin-Lin Hsu in statistics, the paper can
not be completed.

Allow us to discuss more mathematics. Since quantum theory is built on the
Hilbert space, the physicists use the complete sets of the space. Therefore,
the statisticians working on nonparametric approach use the same tool as
physicists. In order to avoid the roughness introduced by the complete sets,
we use the result of Stone-Weirstrass theorem. If it is necessary, then we
will treat the space of continuous functions or measurable functions as
metric space or topological space. Therefore, we use countable dense subset
of the space, the set of Bernstein polynomials.

\textbf{References}

[1]. A. Wald, (1949),'' Note on the Consistency of the Maximum Likelihood
Estimate,'' The Annal of Mathematical Statistics, Vol. 20, No. 4

\noindent
pp. 595-601.

[2]. W. M. Newman and R. F. Sproull, ``Principle of Interactive Computer
Graphics'', McGraw-Hill, New York, (1979), pp. 309-331.

[3]. A. Quarteroni, R, Sacco and F. Saleri, ''Numerical Mathematics'',
Springer, (2000), pp 361-375.

[4]. Tom M. Apostol, `` Calculus '', Vol. 1 John Wiley {\&} Sons, (1967),
pp.374-443.

[5]. Tom M. Apostol, `` Mathematical Analysis'', Addison-Wesely, (1974),
pp.183-247, pp.322.

[6]. R. O. Duda and P. E. Hart, ``Pattern Classification and Scene Analysis
''.John Wiley, (1973), pp. 85-91.

[7] S. Lang, ``linear Algebra'', Springer: 3$^{rd }$ ed., (1976),
pp.125-131.

[8]. D. G. Luenberger, ``Optimization By Vector Space Methods,'' Wiley
(1969), pp.239-265.

[9]. I. J. Good and R. A. Gaskins, Biometrika, Vol. 58, No. 2, (1971), pp.
255-277

[10]. I. J. Good and R. A. Gaskins, Journal of the America Statistical
Association, Vol. 75, No. 369, (1980), pp.42-73.

[11]. M. X. Dong and R. J-B. Wetes, `` Estimating Density Functions: a
Constrained Maximum Likelihood Approach,'' Journal of Nonparametric
statistics, Vol. 12, (2000), pp. 549-595.

[12]. I. A. Ahmad and I. S. Ran, `` Kernel Contrast: A Data-Based Method Of
Choosing Parameters In Nonparametric Density Estimation,'' Journal of
Nonparametric Statistics, Vol. 16(\ref{eq5}), (2004), pp. 671-707.

[13]. B. Hatfield, `` Quantum Field Theory of Particles and Strings,''
(1992), pp.20, pp. 698.

[14]. J. S. Townsend, `` A Modern Approach to Quantum Mechanics'',
McGraw-Hill, (1992), pp.1-24.

\begin{figure}[htbp]
\centerline{\includegraphics[width=5.20in,height=3.67in]{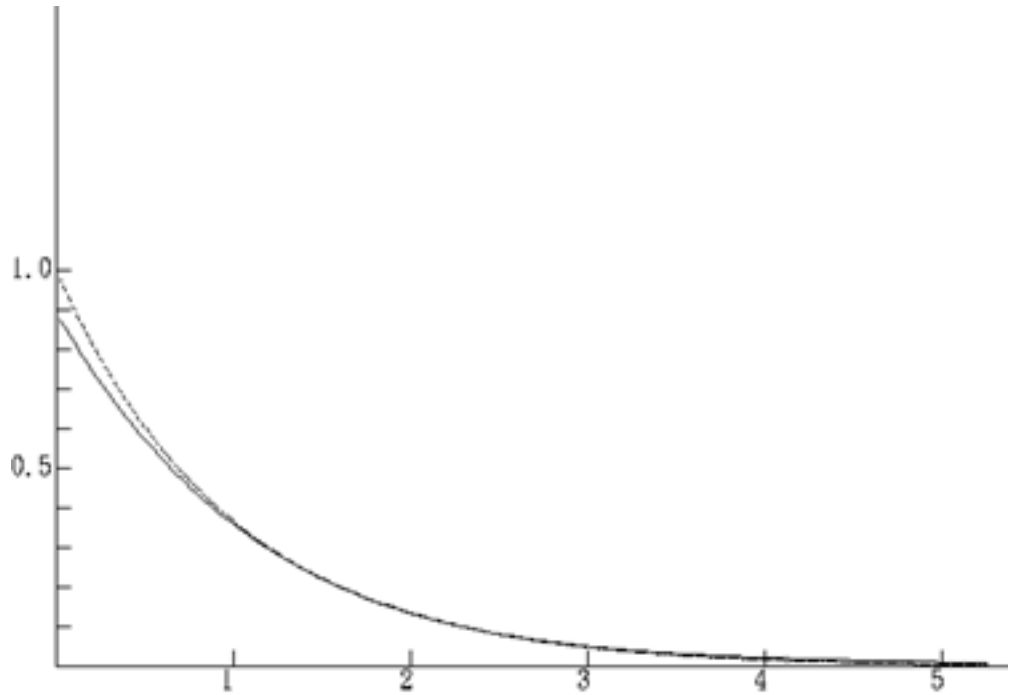}}
\caption{The density function is exp(-x). The sample size is 80.
There are 11 windows of Bezier spline, n=10. }\label{fig1}
\end{figure}
\begin{figure}[htbp]
\centerline{\includegraphics[width=4.92in,height=3.65in]{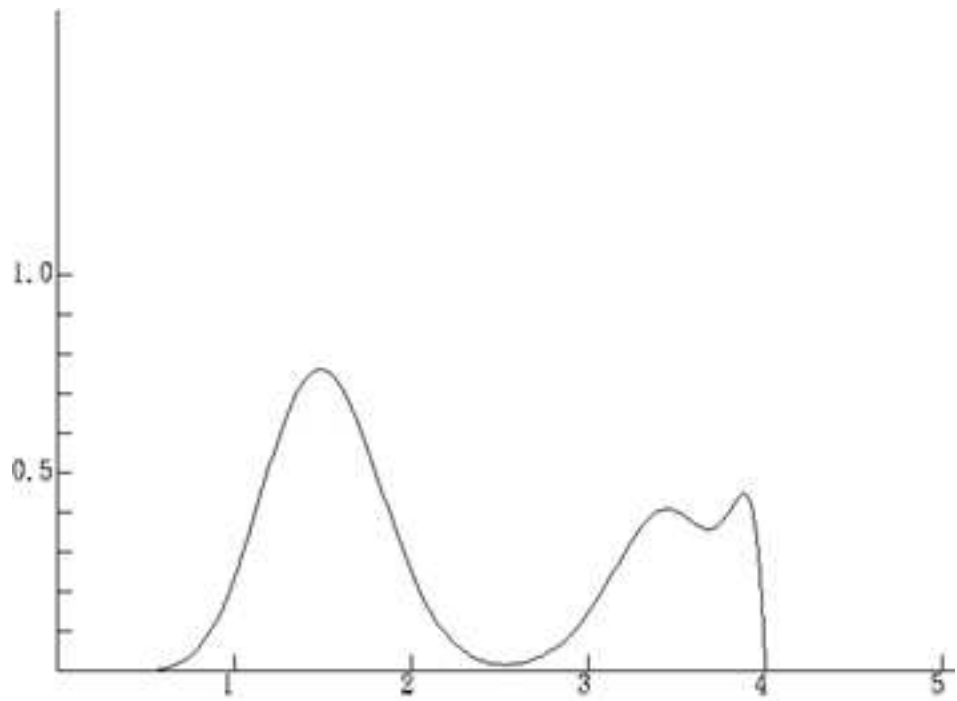}}
\caption{The density function is bimodal. The domain is  [0,4].
f(x)=2/3 when x is in [1,2]; f(x)=1/3 when x is in [3,4]; f(x)=0,
otherwise. The sample size is 180. There are 35 windows of Bezier
spline. } \label{fig2}
\end{figure}
\begin{figure}[htbp]
\centerline{\includegraphics[width=7.17in,height=5.17in]{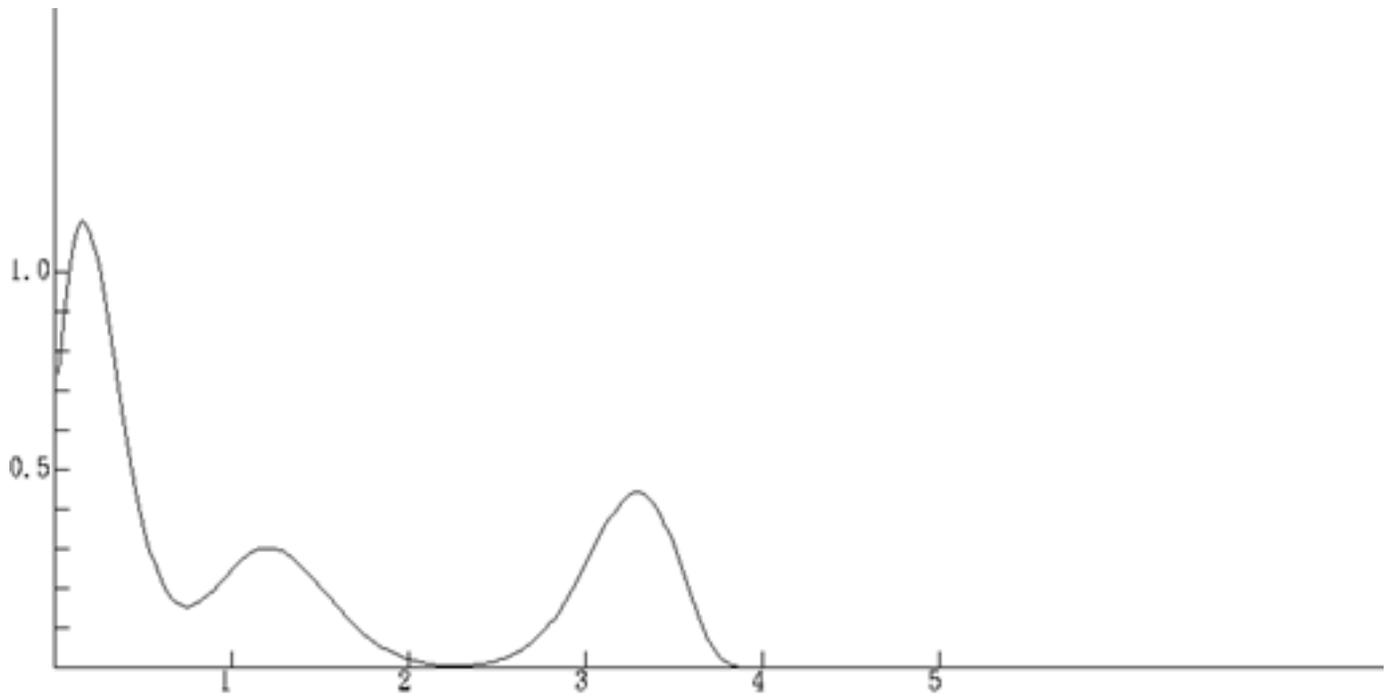}}
\caption{The density function is trimodal.The domain is  [0,4].
f(x)=1 when x is in [0,1/2]; f(x)=1/2 when x is in [1,3/2];
f(x)=1/2 when x is in [3,7/2]; f(x)=0, otherwise. The sample size
is 180. There are 35 windows of Bezier spline. } \label{fig3}
\end{figure}

\end{document}